\documentclass[12pt]{amsart}
\usepackage{amsmath,amssymb}

\newcommand{\Z}{{\mathbf Z}}

\newcommand{\R}{{\mathbf R}}

\newcommand{\supp}{\operatorname{supp}}

\newcommand{\sumstar}{\sum\nolimits^*}

\newcommand{\vol}{\operatorname{vol}}

\newcommand{\meas}{\operatorname{meas}} 

\newcommand{\EX}{\mathbf E} 
\newcommand{\var}{\operatorname{var}}

\numberwithin{equation}{section}

\begin{document}

\newtheorem{thm}{Theorem}[section]

\newtheorem{lem}[thm]{Lemma}
\newtheorem{prop}[thm]{Proposition}
\newtheorem{cor}[thm]{Corollary}

\theoremstyle{definition}
\newtheorem{defn}{Definition}[section]
\newtheorem{conj}{Conjecture}

\theoremstyle{remark}
\newtheorem{remark}{Remark}[section] 
\newtheorem{note}{Note}[section]

\theoremstyle{plain}
\newtheorem{hyp}{H}

\title[spacings of fractional parts of lacunary sequences]
{The distribution of spacings between fractional parts of
lacunary sequences}
\author{Ze\'ev Rudnick and Alexandru Zaharescu}

\address{Raymond and Beverly Sackler School of Mathematical Sciences,
Tel Aviv University, Tel Aviv 69978, Israel 
({\tt rudnick@math.tau.ac.il})}
\address{School of Mathematics, Institute for Advanced Study, Olden Lane, 
Princeton, NJ 08540 ({\tt zaharesc@math.ias.edu})}

\date{December 12, 1999}
\thanks{Supported in part by a grant from 
the U.S.-Israel bi-national Science Foundation} 
\maketitle

\section{Introduction}

\subsection{}
A {\em lacunary sequence} is a sequence of integers $a(x)$, $x=1,2,\dots$
which satisfies the ``gap condition''
$$
\liminf \frac{a(x+1)}{a(x)}>1 \;.
$$
A primary example is to take an integer $g\geq 2$ and set $a(x)=g^x$. 

As is true for any increasing sequence of integers, for almost every $\alpha$
the fractional parts $\alpha a(x)$ are uniformly distributed modulo 1.
Moreover, for lacunary sequences, it has long been known that
the fractional parts of $\alpha a(x)$ have strong randomness properties.
For instance,  the exponential sums
$\frac 1{\sqrt{N}}\sum_{x\leq N} \cos(2\pi  \alpha a(x))$
have a Gaussian value distribution as $N\to\infty$
(see the survey in \cite{Kac}).

In this paper, we show that lacunary sequences have  additional
features in common with those of random sequences, which is the
asymptotic distribution of {\em spacings} between elements of the
sequence: 
Given a  sequence $\{\theta_n\} \subset [0,1)$, the
nearest-neighbor spacing distribution is defined by ordering the first
$N$ elements of the sequence:
$\theta_{1,N}\leq \theta_{2,N}\leq \dots \leq \theta_{N,N}$, and then defining the
normalized spacings to be
$$
\delta_n^{(N)}:= N(\theta_{n+1,N} - \theta_{n,N}) \;.
$$
The asymptotic distribution function of $\{\delta_n^{(N)} \}_{n=1}^N$
is level spacing distribution $P_1(s)$, that is for each interval
$[a,b]$ we require that
$$
\lim_N \frac 1N\#\{n< N: \delta_n^{(N)}\in [a,b] \} =
\int_a^b P_1(s)ds \;.
$$

The statistical model we have in mind is the
``Poisson model'', of a sequence generated by uncorrelated levels
(i.i.d.'s). In that case $P(s)=e^{-s}$. Moreover in that model one
knows the behavior of all other local spacing statistics, such as for
instance:
\begin{enumerate}
\item
Instead of spacings between nearest neighbors, one can consider
spacings between next-to-nearest neighbors or more generally for any
fixed $a\geq 1$, set
$$
\delta_{a,n}^{(N)}:=N(\theta^{(N)}_{n+a}-\theta^{(N)}_n) \;.
$$
Let $P_a(s)$ be the limiting distribution function of
$\{\delta_{a,n}^{(N)}\}$ as $N\to \infty$. In the Poisson model,
$$P_a(s) = \frac{s^{a-1}}{(a-1)!}e^{-s} \;. $$
\item
For fixed $r\geq 1$ consider the joint distribution of the nearest
neighbor spacings
$(\delta_{n}^{(N)}, \delta_{n+1}^{(N)}, \dots ,\delta_{n+r-1}^{(N)})$.
In the Poisson model, these are independent and so the distribution
function is $\prod_{i=1}^{r} e^{-s_i}$.
\item
For fixed $\lambda > 0$, consider the probability of finding exactly 
$k$ elements of the sequence $\{\theta_n: n\leq N\}$ 
in a randomly chosen interval of length $\lambda/N$. 
In the Poisson model, this probability is $e^{-\lambda}\frac{\lambda^k}{k!}$.  
\end{enumerate}

\subsection{Results} 
The principal result of our paper asserts that
\begin{thm}\label{ae for P}
Let $a(x)$ be a lacunary sequence. Then for almost all
$\alpha$,  the fractional parts of the sequence $\{\alpha a(x)\}$
has all its local spacing measures as those of the Poisson model.
\end{thm}

As is well known
, all local spacing measures are determined by
the correlation functions, which measure the distribution of spacings
between tuples of elements, not necessarily neighboring.
To define the $k$-level correlation function,
for $x=(x_1,\dots,x_k)$, denote by $\Delta(x)$ the  difference vector
$$
\Delta(x) = (a(x_1)-a(x_2),\dots, a(x_{k-1})-a(x_k) ) \;.
$$
Take a smooth, compactly supported function $f\in
C_c^\infty(\R^{k-1})$, and set
$$F_N(y):= \sum_{m\in \Z^{k-1}} f(N(m+y)) \;.
$$
We then define the $k$-level correlation sum associated to this data by
\begin{equation}\label{definition of R_k}
R_k(f,N)(\alpha):=\frac 1N\sumstar_{x_i\leq N} F_N(\alpha \Delta(x))
\end{equation}
where $\sumstar$ means the sum over all vectors with {\em distinct}
components: $x_i\neq x_j$ if $i\neq j$.
Our main result is:
\begin{thm}\label{ae for R_k}
There is a set of  $\alpha$ of full measure so that for
all $k\geq 2$ and all test functions $f\in C_c^\infty(\R^{k-1})$, the
$k$-level correlation sums $R_k(f,N)(\alpha)$ converge to
$\int f(x)dx$.
\end{thm}
By standard results
, this implies Theorem~\ref{ae for P}. 
The case of  pair correlation  ($k=2$) was done in \cite{RZ}.

\subsection{Comparison with polynomial sequences} 
Much of the work  done previously on spacings of fractional parts was
for polynomial sequences, such as $a(x)=x^2$ 
\cite{BZ, RS, RSZ}, see also \cite{MS, Zelditch}.  
Rudnick and Sarnak \cite{RS} proved the analogue of
Theorem~\ref{ae for R_k} for the pair correlation function ($k=2$). However, 
the method used both in \cite{RS} and here, 
which proves almost-everywhere convergence by going through
convergence in $L^2$, already fails in the case of $a(x)=x^2$  
at the level of triple correlation, 
because the variance diverges as $N\to \infty$.

The  reason for the difference  between these two cases 
can be understood by examining the number of solutions of the equation
\begin{multline}\label{tripeq}
n_1( a(x_1)-a(x_2) ) + n_2( (a(x_2)-a(x_3) ) = \\
n'_1( a(x'_1)-a(x'_2 )) + n'_2( a(x'_2)-a(x'_3) ) 
\end{multline}
in variables bounded by $N$, and $n,n'\neq 0$. 
For $a(x)=x^2$ the
number of solutions of \eqref{tripeq} is $\gg N^7$. 
This is consistent with the heuristic that zero is a typical value 
of the difference of the two sides of the equation, and for $a(x)$
growing as slowly as $x^2$ the size
of this difference  is at most $O(N^3)$ while the number
of variables is $10$. Thus the typical difference should occur about
$N^7$ times.  
As is explained in \cite{RS}, 
this effect causes the variance of $R_3(f,N)$ to
blow up like $N$. 
A similar effect will cause the blow-up of the variance of 
high correlations  for any polynomially increasing sequence. 
The non-Gaussian distribution of the ``theta sums'' 
$\frac 1{\sqrt{N}}\sum_{x\leq N}\exp(2\pi i \alpha x^2)$ is related to  
this kind of clustering effect \cite{Jurkat, Marklof}.

In contrast, for lacunary sequences we will show in
section~\ref{sec:counting} that the number of solutions of \eqref{tripeq} 
is $O( (N\log N)^5 )$, which is not much more than the number of
``diagonal'' solutions $x=x'$, $n=n'$.

\subsection{Plan of the paper} 
We begin in section~\ref{sec:counting} with a key counting argument: 
We consider  the number of solutions of an equation 
\begin{multline}\label{sysvar}
m_1(a(x_1)-a(x_2))+\cdots +m_{k-1}(a(x_{k-1})-a(x_k))\\
=m'_1(a(x'_1)-a(x'_2))+\cdots +m'_{k-1}(a(x'_{k-1})-a(x'_k)).
\end{multline}
in integers  $0\neq m,m'\in [-N,N]^{k-1}$, $x,x'\in [1,N]^k$,  
$x_1,\cdots ,x_k$ distinct, $x_1',\cdots ,x_k'$ distinct.  
In Lemma~\ref{lem z4} 
we show  that the number of such solutions is
$O(N^{2k-1}\log^{2k-1}N)$. This is comparable to the 
number of ``diagonal'' solutions, which is of order $N^{2k-1}$. 
For {\em fixed} coefficients $m$, $m'$, the
diagonal solutions are indeed responsible for the bulk of the
solutions, see e.g. \cite{Israilov}. 

We then show in section~\ref{sec:mean}  that the mean of
$R_k(f,N)$ is asymptotic to $\int f$, and in
section~\ref{sec:variance} we show that the variance decays with $N$: 
$\var(R_k(f,N))\ll N^{-1+\epsilon}$, for all $\epsilon>0$. These are
done by a reduction to the study of solutions of \eqref{sysvar}. 

In section~\ref{sec:ae} we show almost-everywhere convergence, after first  
investigating in section~\ref{sec:small} the frequency of occurrence of
fractional parts of $\alpha a(x)$ in short (of size $1/N$) intervals.

\section{A counting lemma}\label{sec:counting}

Let $a(x)$ be a lacunary sequence, that is there is some  $c>1$ so that 
$$a(n+1)>ca(n)$$ 
for all $n$ sufficiently large.
We wish to estimate the number of solutions of an equation such as
\eqref{sysvar}. We will do so in  Lemma~\ref{lem z4}, after some
preliminaries. 

\begin{lem} \label{lem z1}
Let $s\ge 1$, $C > 0$ and let $A_1>A_2>\cdots >A_s $ be positive integers.
Then for any $b\in \bf Z $ and $N\ge 1$ the number of vectors 
$\vec y =(y_1,\cdots ,y_s)\in \bf Z^s$ with $|y_1|,\cdots ,|y_s|\le N$
such that
\begin{equation}\label{sandwich}
|y_1A_1+\cdots +y_sA_s+b|\le CA_1
\end{equation}
is $O_{s,C}(N^{s-1})$.
\end{lem}
\begin{proof}  We need to count the number of integer points $\vec y$ inside 
the region $\Omega \subset \bf R^s$ which consists of the points in the cube 
$[-N,N]^s$ which lie between the hyper-planes

\begin{equation}\label{hyperplanes}
\begin{split}
y_1A_1+\cdots +y_sA_s+b= CA_1\\
y_1A_1+\cdots +y_sA_s+b= -CA_1\;.
\end{split}
\end{equation}
Note that the region $\Omega $ is convex and contained in a ball
around the origin of radius $\ll _s N$.
By the Lipschitz principle (see \cite{Davenport}) we know that
\begin{equation}\label{lips}
\#(\Omega \cap {\bf Z}^s)=\vol(\Omega )+O_s(N^{s-1}).
\end{equation}

The distance between the above hyper-planes is  
$$
\frac {2CA_1}{\sqrt {A_1^2+\cdots +A_s^2}}\le 2C 
$$
thus $\Omega $ is contained in a cylinder of height $2C$ whose base is an
$(s-1)-$ dimensional ball of radius $\ll_s N$. 
Therefore $\vol(\Omega )=O_{s,C}(N^{s-1})$
which together with \eqref{lips} gives the lemma.
\end{proof}

\begin{lem}  \label{lem z2}
Let $s\ge 2 $ and $z_1>\cdots >z_s $ be positive integers.
Then for any $b,d\in \bf Z$ and any $N\ge 1$ the number of vectors 
$\vec y=(y_1,\cdots ,y_s)\in \bf Z^s$ with 
$|y_1|,\cdots ,|y_s|\le N$ for which
\begin{equation}\label{linear}
\begin{split}
|y_1a(z_1)+\cdots +y_s a(z_s)+b|&\le C a(z_1) \\
y_1+\cdots +y_s+d&=0
\end{split}
\end{equation}
holds true is $O_{s,c}(N^{s-2})$.
\end{lem}
\begin{proof}  
We first remark that since $s\ge 2$ and the $z's$ are distinct,
 the hyper-planes \eqref{hyperplanes} with $A_1,\cdots,A_s$ replaced
replaced by $a(z_1),\cdots,a(z_s)$
are not parallel to the hyper-plane  given by the equation \eqref{linear}.
Moreover, the fact that our sequence is lacunary 
insures that the angle between these 
hyper-planes is not small. 
Thus when we solve for $y_s$ in \eqref{linear} and input the result 
in \eqref{sandwich} we get an inequality in $s-1$ variables:
\begin{equation}\label{ineq in s-1}
|y_1(a(z_1)-a(z_s))+\cdots +y_{s-1}(a(z_{s-1})-a(z_s))+b-da(z_s)|\le a(z_1)
\end{equation}
in which the RHS is bounded by the largest of the coefficients which appear
in the LHS: 
$$
a(z_1)-a(z_s)\ge (1-\frac 1c)a(z_1)\;.
$$
Then Lemma~\ref{lem z1} applies to \eqref{ineq in s-1}, with
$A_j=a(z_j)-a(z_s)$ for $1 \le j \le s-1$ and $C=\big (1- \frac 1c\big )^{-1}$,
and we find that the number of vectors $\vec y$
having the required properties is $O_{s,c}(N^{s-2})$ as stated.
\end{proof}

We now come to our main counting lemma.
\begin{lem} \label{lem z3}
Let $r\ge 1 $ be an integer. For any $N\ge 1 $ the number
of solutions $(y_1,\cdots ,y_r,z_1,\cdots ,z_r)$ to the system:
\begin{equation}\label{system}
\begin{split}
y_1a(z_1)+\cdots +y_ra(z_r)&=0\\
y_1+\cdots +y_r&=0
\end{split}
\end{equation}
in integers  $y_1,\cdots,y_r$,  
$$
(y_1,\cdots y_r)\ne (0,\cdots ,0)
$$
$$
z_1,\cdots ,z_r \geq 1 \quad\text{ distinct }
$$
$$
|y_1|,\cdots ,|y_r|,|z_1|,\cdots ,|z_r|\le N
$$
is $O_{r,c}(N^{r-1}\log^{r-1}N)$.
\end{lem}
\begin{proof}  
Our proof is by induction on $r$. The case $r=1$ is clear,
the number of solutions
in this case being zero. Let us assume that the statement holds true for $r-1$ 
and prove it for $r$. Let $(y_1,\cdots ,y_r,z_1,\cdots ,z_r)$ be a solution to the 
system \eqref{system}. If there exists $j\in \{1,\cdots ,r\}$ such
that $y_j=0$ then  
$(y_1,\cdots ,y_{j-1},y_{j+1},\cdots ,y_r,z_1,\cdots ,z_{j-1},z_{j+1},\cdots ,z_r)$
will be a solution for the same system with $r$ replaced  by $r-1$.
By the induction assumption the number of solutions of this system is 
$O_{r,c}(N^{r-2}\log^{r-2}N)$. For each such solution, $z_j$ is free to take 
values $\le N$. Therefore the number of solutions to the system
\eqref{system}  for which
at least one of $y_1,\cdots ,y_r$ vanishes is $O_{r,c}(N^{r-1}\log^{r-2}N)$.
We now count the solutions to \eqref{system} with $y_j\ne 0$ for all $j$.
There are $r!$ possible orders for the $z's$. Let 
us count the solutions for which $z_1>\cdots >z_r$. Given  such a solution
$(y_1,\cdots ,y_r,z_1,\cdots ,z_r)$ we consider the partition of the set
$\{1,\cdots ,r\}$ as a disjoint union of sets $B_1,\cdots ,B_l$ defined as 
follows. $B_1$ consists of those $j\in \{1,\cdots ,r\} $ for which $z_j\ge z_1-
\frac {2\log N}{\log c}$.
If $j_2$ is the smallest index not contained in $B_1$ then we put in $B_2$ all those 
$j\in \{j_2,\cdots ,r\}$ for which $ z_j\ge z_{j_2}-
\frac {2\log N}{\log c}$, and so on . In the end, if $1=j_1<j_2<\cdots <j_l$ are the 
smallest indices contained in $B_1,B_2,\cdots ,B_l$
respectively, then we have:
\begin{equation}\label{windows}
z_{j_2}<z_{j_1}-\frac {2\log N}{\log c}\le z_{j_2-1},\cdots ,z_{j_l}<
z_{j_{l-1}}-\frac {2\log N}{\log c}\le z_{j_{l}-1}.
\end{equation}

The number of partitions as above is bounded in terms of $r$. 
Let us count the number of solutions $(y_1,\cdots ,y_r,z_1,\cdots ,z_r)$ 
which correspond to a given partition $B_1,\cdots ,B_l$.
We distinguish two cases: $\# B_l\ge 2$ and $\# B_l=1$.

Let us first treat the case $\# B_l\ge 2$. If we fix 
$z_{j_1},z_{j_2},\cdots z_{j_l}$ then from \eqref{windows} it follows
that each  
of the remaining $z's$ can take at most $[\frac {2\log N}{\log c}]$ values.
Hence the number of vectors $\vec z =(z_1,\cdots ,z_r)$ satisfying
\eqref{windows} is $O_{r,c}(N^l\log^{r-l}N)$.
Thus we are done with the case $\# B_l\ge 2$ if we show that for any vector 
$\vec z$ as above the number of solutions $\vec y=(y_1,\cdots ,y_r)$ is 
$O_{r,c}(N^{r-l-1})$. Fix some such $\vec z$ and note that by
\eqref{windows} one has:
\begin{equation}\label{bddd}
a(z_{j_2})<\frac {a(z_1)}{N^2},\cdots ,a(z_{j_l})<\frac {a(z_{j_{l-1}})}
{N^2}.
\end{equation}
Let us take a solution $\vec y$ and look at its first $j_2-1$
components.
These are nonzero integer numbers in the interval $[-N,N]$
satisfying the inequality:
$$
|y_1a(z_1)+\cdots +y_{j_2-1}a(z_{j_2-1})|=|\sum _{j\ge j_2}
y_ja(z_j)|<rNa(z_{j_2})<a(z_1).
$$
Here we may apply Lemma~\ref{lem z1} with $s=j_2-1$,
$b=0$  and $A_1,\cdots, A_s$ replaced by $a(z_1),\cdots,a(z_s)$
to conclude that the vector 
$(y_1,\cdots ,y_{j_2-1})$ can only take $O_r(N^{j_2-2})$ values.
Let us fix $(y_1,\cdots ,y_{j_2-1})$ and count the number of 
solutions $\vec y $ whose first $j_2-1$ components are 
$y_1,\cdots ,y_{j_2-1}$. We are now interested in those 
components $y_j$ of $\vec y$ for which  $j\in B_2$.
Write $b=y_1a(z_1)+\cdots +y_{j_2-1}a(z_{j_2-1})$ and use
\eqref{bddd} to deduce that for any solution $\vec y$, its components
$y_j$ with $j\in B_2$ satisfy the inequality:
$$
|b+\sum _{j\in B_2}y_ja(z_j)|=|\sum _{j\ge j_3}y_ja(z_j)|<a(z_{j_2}).
$$
By Lemma~\ref{lem z1} we know that as $\vec y$ varies, the vector formed with the 
components $y_j$ of $\vec y$ for $j\in B_2$ can only take 
$O_r(N^{\# B_2-1})$ values.
We now repeat the above reasoning until we get to the last set of 
components of $\vec y$, namely the $y_j's$ with $j\in B_l$.
The components $y_j$ with $j<j_l$ being fixed, write 
$b=\sum _{1\le j<j_l}y_ja(z_j),d=\sum _{1\le j<j_l}y_j$
and then apply Lemma~\ref{lem z2} (here one uses the assumption that $\# B_l\ge 2$).
It follows that the vector formed with the components $y_j,j\in B_l$
of $\vec y$ can take $O_{r,c}(N^{\# B_l-2})$ values only.
The number of solutions $\vec y $ for a fixed $\vec z$ as above is then 
$\ll _{r,c}N^{(\# B_1-1)+\cdots +(\# B_{l-1}-1)+(\# B_l-2)}=N^{r-l-1}$,
which completes the proof in case $\# B_l\ge 2$.

Assume now that $\# B_l=1$. Then $j_l=r$. In this case we fix 
$z_1,\cdots ,z_{r-1}$ only . This can be done in $O_{r,c}(N^{l-1}\log^{r-l}N)$
ways. For $z_1,\cdots ,z_{r-1}$ fixed we apply Lemma~\ref{lem z1} repeatedly
to conclude that as the vector $(y_1,\cdots ,y_r,z_r)$ varies in the set 
of solutions, the vector $(y_1,\cdots ,y_{r-1})$ can take 
$O_{r,c}(N^{(\# B_1-1)+\cdots +(\# B_{l-1}-1)})=O_{r,c}(N^{r-l})$ values only.
Now for $y_1,\cdots ,y_{r-1},z_1,\cdots ,z_{r-1}$ fixed, $y_r$ and $z_r$  are
uniquely determined from the last two relations in \eqref{system}
(here one uses the fact that 
$y_r\ne 0$). 
Thus the number of solutions $(y_1,\cdots ,y_r,z_1,\cdots ,z_r)$ is 
$O_{r,c}(N^{r-1}\log^{r-l}N)$ in case $\# B_l=1$ as well, and the lemma is
proved. 
\end{proof}

We intend to use the above counting lemma to bound the number of solutions 
of the following equation:
\begin{multline}\label{sys2}
m_1(a(n_1)-a(n_2))+\cdots +m_{k-1}(a(n_{k-1})-a(n_k))\\
=m'_1(a(n'_1)-a(n'_2))+\cdots +m'_{k-1}(a(n'_{k-1})-a(n'_k)).
\end{multline}
in variables $m_1,\cdots ,m_{k-1},m_1',\cdots ,m_{k-1}'\in {\Z}$, 
$n_1,\cdots ,n_k,n_1',\cdots ,n_k'\in {\bf N}$, $n_1,\cdots ,n_k$
distinct, $n_1',\cdots ,n_k'$ distinct, 
$$
(m_1,\cdots ,m_{k-1},m_1',
\cdots ,m_{k-1}')\ne (0,\cdots ,0)
$$
and all variables of absolute value at most $N$.

The result we obtain is the following:
\begin{lem} \label{lem z4}
Let $k\ge 2, k\in {\bf Z}$. For any $N\ge 1$ the number of solutions  to 
the system \eqref{sys2} is $O_{k,c}(N^{2k-1}\log^{2k-1}N)$.
\end{lem}

In order to simplify the combinatorics involved in the derivation of
Lemma~\ref{lem z4} 
from Lemma~\ref{lem z3} we first establish a more general form of
Lemma~\ref{lem z3}.
Let $r\ge 1$ and consider the system:
\begin{equation}\label{sys3}
\begin{split}
y_1a(z_1)+\cdots +y_ra(z_r)&=0\\
y_1+\cdots +y_r&=0
\end{split}
\end{equation}
$$
y_1,\cdots ,y_r\in {\bf Z}, \quad z_1,\cdots ,z_r\in {\bf N}
$$
$$
|y_1|,\cdots ,|y_r|,|z_1|,\cdots ,|z_r|\le N
$$

Let $(\vec y,\vec z )=(y_1,\cdots ,y_r,z_1,\cdots ,z_r)$ be a solution of \eqref{sys3}.
For any $i\in \{1,\cdots ,r\} $ denote $A(i)=\{1\le j\le r: z_j=z_i\}$.
We say that the solution $(\vec y,\vec z)$ is {\em degenerate} 
provided we have 
\begin{equation}\label{Degenerate}
\sum_{j\in A(i)}y_j=0
\end{equation}
for all $i\in \{1,\cdots ,r\}$. Otherwise we say that $(\vec y,\vec z)$
is {\em non-degenerate}. We have the following :

\begin{lem} \label{lem z5}
Let $r\ge 1$. Then for any $N\ge 1$ the number of non-degenerate 
solutions to \eqref{sys3} is $O_{r,c}(N^{r-1}\log^{r-1}N)$.
\end{lem}
\begin{proof}  Each solution $(\vec y,\vec z)$ to \eqref{sys3} produces a
partition of the set  
$\{1,\cdots ,r\}$ as a disjoint union of subsets $A_1,\cdots ,A_l$,
where $A_1,\cdots ,A_l$
are the above sets $A(1),\cdots ,A(r)$ without repetitions.
Let us count the number of non-degenerate solutions to \eqref{sys3}
which correspond  to 
a given partition $A_1,\cdots ,A_l$ of the set $\{1,\cdots ,r\}$. For $s=1,2,
\cdots ,l$ denote $u_s=\sum _{j\in A_s}y_j$, $v_s=z_j$ for $j\in A_s$, then write
$\vec u=(u_1,\cdots ,u_s),\vec v=(v_1,\cdots ,v_s)$.
If $(\vec y,\vec z)$ is a non-degenerate solution to \eqref{sys3} then
not all the numbers  
$u_1,\cdots ,u_s$ vanish. One sees that for any such $(\vec y,\vec z )$
the pair $(\vec u,\vec v)$ is a solution of the system:
\begin{equation}\label{sys4}
\begin{split}
u_1a(z_1)+\cdots +u_la(z_l)&=0 \\
u_1+\cdots +u_l&=0
\end{split}
\end{equation}
in integers $u_1,\cdots ,u_l\in {\bf Z}$, $\vec u\ne \vec 0$, 
$v_1,\cdots ,v_l\in N$  distinct, 
$$
|u_1|,\cdots ,|u_l|,|v_1|,\cdots ,|v_l|\le N
$$

By Lemma~\ref{lem z3}  we know that the number of solutions of the system
\eqref{sys4} 
is $O_{l,c}(N^{l-1}\log^{l-1}N)$. Now fix a solution $(\vec u, \vec v)$
and count the number of non-degenerate solutions $(\vec y,\vec z)$
to \eqref{sys3} which correspond to the above partition $A_1,\cdots ,A_l$
and which produce the vector $(\vec u, \vec v)$. Clearly $\vec z$ is 
uniquely determined since $z_j=v_s$ for any $s$ and any $j\in A_s$.
Moreover, for any $s$ the number of solutions $y_j,j\in A_s$ of the 
equation $\sum _{j\in A_s}y_j=u_s$ is $O_r(N^{\# A_s-1})$.
Hence the number of solutions $(\vec y,\vec z)$ which correspond to a 
given pair $(\vec u, \vec v)$ is 
$O_r(N^{(\# A_1-1)+\cdots +(\# A_s-1)})=O_r(N^{r-l})$
and so the total number of non-degenerate solutions to \eqref{sys3} is 
$O_r(N^{r-1}\log ^{r-1}N)$, which completes the proof of
Lemma~\ref{lem z5} .
\end{proof}

\noindent{\bf Proof of Lemma~\ref{lem z4}:}  
Denote $r=2k, z_1=n_1,\cdots ,
z_k=n_k, z_{k+1}=n'_1,\cdots ,z_r=n'_k, y_1=m_1, y_2=m_2-m_1,
\cdots ,y_{k-1}=m_{k-1}-m_{k-2}, y_k=-m_{k-1}, y_{k+1}=-m'_1,
y_{k+2}=m'_1-m'_2,\cdots ,y_{2k-1}=m'_{k-2}-m'_{k-1}$
and $y_k=m'_{k-1}$. Then any solution $(\vec m,\vec n,\vec m',\vec n')$
of \eqref{sys2} produces a solution $(\vec y,\vec z)$ of \eqref{sys3}
(with $N$ replaced  
by $2N$) which satisfies the additional properties:
\begin{equation}\label{additional}
y_1+\cdots +y_k=0
\end{equation}
with $(y_1,\cdots ,y_r)\ne (0,\cdots ,0)$, $z_1,\cdots ,z_k$ distinct,
$z_{k+1},\cdots z_r$  distinct, 
and each such $(\vec y,\vec z)$ uniquely determines the tuple 
 $(\vec m,\vec n,\vec m',\vec n')$.

Thus we are done if we show that the number of solutions to \eqref{sys3} which 
satisfy the additional requirements \eqref{additional} is
$O_r(N^{r-1}\log ^{r-1}N)$. 
Lemma~\ref{lem z5} takes care of the non-degenerate solutions to \eqref{sys3}
 so it remains  
to count the number of degenerate solutions to \eqref{sys3} which satisfy 
\eqref{additional}.

Let $(\vec y,\vec z)$ be such a solution. If $z_1,\cdots ,z_r$ are 
distinct then by the degeneracy conditions \eqref{Degenerate} it follows that 
$y_1=y_2=\cdots =y_r=0$ which contradicts \eqref{additional}. Thus some $z_j$
with $1\le j\le k$ will have to equal some $z_j$ with $ k+1\le j\le 2k$.
Let $s$ be the number of indices $j\in \{1,\cdots ,k\}$ for which there exists 
$i\in \{k+1,\cdots ,2k\}$ such that $z_j=z_i$. Both \eqref{sys3} and
\eqref{additional}   
are symmetric in $z_1,\cdots ,z_k$ and separately in $z_{k+1},\cdots ,
z_{2k}$ and the same holds true for $y_1,\cdots ,y_k$ respectively 
$y_{k+1},\cdots ,y_{2k}$. After making a permutation of variables if 
necessary, we may assume that $z_j=z_{j+k}$ for $1\le j\le s$.
Then the sets $A_1,\cdots ,A_l$  look like this : $A_1=\{1,k+1\},
A_2=\{2,k+2\},\cdots , A_s=\{s,k+s\}, A_{s+1}=\{s+1\},\cdots ,
A_k=\{k\}, A_{k+1}=\{k+s+1\},\cdots ,A_l=\{2k\}$, where $l=2k-s$.
The degeneracy relations \eqref{Degenerate} become:
\begin{equation}\label{newdeg}
\begin{cases}
y_j+y_{j+k}=0,& 1\le j\le s, \\
y_j=0,&  s+1\le j\le k  \quad\text{ or }  k+s+1\le j\le 2k \;. 
\end{cases}
\end{equation}
Now, given an $s\in \{1,\cdots ,k\}$ and the above partition 
$A_1,\cdots ,A_l$, the number of degenerate solutions 
$(\vec y,\vec z)$  which correspond to this partition are counted as
follows.
On one hand each of the $l$ distinct $z's$ can assume at most $N$ values ,
so $\vec z$ takes at most $N^l=N^{2k-s}$ values. On the other hand,
each of the variables $y_j$ (if there are any) with $2\le j\le s $ 
assumes at most $2N+1$ values and for each such choice of the vector
$(y_2,\cdots ,y_s)$ the variables $y_{k+2},\dots ,y_{k+s},y_{s+1},
\dots ,y_k,y_{k+s+1},\dots ,y_{2k}$ are determined by \eqref{newdeg},
then $y_1$ is determined by \eqref{additional} and the remaining variable 
$y_{k+1}$ is determined by \eqref{newdeg}.
Hence $\vec y$  takes at most $ (2N+1)^{s-1}$ values and the number of 
degenerate  solutions $(\vec y,\vec z)$ is $O_r(N^{2k-1})$,which completes
the proof of Lemma~\ref{lem z4}. \qed


\section{The average value of $R_k(f,N)$} \label{sec:mean}
\subsection{Poisson sum} 
Recall that for $f\in C_c^\infty(\R^{k-1})$, $y\in \R^{k-1}$, we set 
$$
F_N(y) = \sum_{m\in \Z^{k-1}} f(N(y+m))\;.
$$
By Poisson summation,
\begin{equation}\label{poissum}
F_N(y) = \frac 1{N^{k-1}} \sum_{n\in \Z^{k-1}} 
\widehat f(\frac nN)e(n\cdot y) \;.
\end{equation}
By inserting \eqref{poissum} into the definition 
\eqref{definition of R_k} of $R_k(f,N)$  we find: 
\begin{equation}
R_k(f,N)(\alpha) = \frac 1{N^k}\sum_{n\in \Z^{k-1}} 
\widehat f(\frac nN) \sumstar_{x_i\leq N} e(\alpha n\cdot \Delta(x)) \;.
\end{equation} 

Since $R_k(f,N)(\alpha)$  is periodic in $\alpha$, we may expand it in
a Fourier series
\begin{equation} \label{R_k in terms of b_l}
R_k(f,N)(\alpha) = \frac 1{N^k} \sum_{l\in \Z} b(l,N)e(l\alpha)
\end{equation}
where 
$$
b(l,N) = \sum_{n\in \Z^{k-1}} \sumstar_{\substack{x_i\leq N\\n\cdot
\Delta(x) = l}}    \widehat f(\frac nN) \;.
$$

\subsection{The mean of $R_k(f,N)$} 
 From \eqref{R_k in terms of b_l} we can immediately compute the mean 
of $R_k(f,N)$ as  
$$
\EX(R_k(f,N)) = \int_{0}^1 R_k(f,N)(\alpha)  = \frac{b(0,N)}{N^k}
$$

\begin{lem} \label{lem:average}
Assume $a(x)$ is a lacunary sequence. 
Then $\forall \epsilon>0$,  
$$
\EX(R_k(f,N))= \frac{b(0,N)}{N^k} =
\widehat f(0)+O_{f,\epsilon}(\frac 1{N^{1-\epsilon}})
$$
\end{lem}
\begin{proof}  
We write 
\begin{equation*}
\begin{split}
b(0,N)&=\widehat f(0)\#\{x_i\leq N: \text{distinct}\} + 
\tilde b(N)\\ 
&= \widehat f(0)N^k\left(1+O(\frac 1N)\right) + \tilde b(N)
\end{split}
\end{equation*}
where 
\begin{equation}\label{tide b}
\tilde b(N) = \sum_{n\neq 0} \sumstar_{\substack{x_i\leq N\\n\cdot
\Delta(x) = 0}}    \widehat f(\frac nN) 
\end{equation}
we will show that $\tilde b(N)\ll N^{k-1+\epsilon}$ and thus prove our
lemma. 

Fix $\epsilon>0$, and let $\delta=\epsilon/2(k-1)$, 
$R\geq (100+k)/\delta+ k $. 
Since $f\in C_c^\infty(\R^{k-1})$,
$|\widehat f(x)| \ll |x|^{-R}$ for large $|x|$. 
Now divide the range of summation in \eqref{tide b} into 
$0<|n|\leq N^{1+\delta}$ and $|n|>N^{1+\delta}$:
$$
\tilde b(N) \ll_f \sum_{0<|n|\leq N^{1+\delta}}
\sumstar_{\substack{x_i\leq N\\n\cdot \Delta(x) = 0}} 1 +
\sumstar_{x_i\leq N}\sum_{|n|>N^{1+\delta}}|\frac nN|^{-R}
$$
The second sum is bounded by 
$$
N^{k+R}\sum_{|n|>N^{1+\delta}}\frac 1{|n|^R} \ll
N^{k+R-(1+\delta)(R-k)} \ll N^{k-100}
$$
by our choice of $\delta$ and $R$. 

As for the first sum, it is bounded by the number of 
$x=(x_1,\dots ,x_k)$ with distinct $x_i\leq N^{1+\delta}$, 
and $n\in \Z^{k-1}$ with $0<|n|\leq N^{1+\delta}$ such that 
$n\cdot \Delta(x) = 0$.  
By  Lemma~\ref{lem z3}, this number is 
$\ll (N^{1+\delta}\log (N^{1+\delta}) )^{k-1}\ll N^{k-1+\epsilon}$. 
Thus we find  that $\tilde b(N)\ll N^{k-1+\epsilon}$ as required.
\end{proof}


\section{Estimating the variance} \label{sec:variance}

\begin{prop}\label{prop:variance}
The variance of $R_k(f,N)$ satisfies 
$$
\var(R_k(f,N)):= \int_0^1 \left| R_k(f,N)(\alpha) -\EX(R_k(f,N) \right|^2
d\alpha \ll_\epsilon \frac 1{N^{1-\epsilon}}
$$
for all $\epsilon>0$. 
\end{prop}
\begin{proof}
By \eqref{R_k in terms of b_l} we have 
\begin{equation}\label{var formula}
\begin{split}
\var(R_k(f,N))  &= \EX(\left| R_k(f,N)-b(0,N)\right|^2) \\
&= \frac 1{N^{2k}}\sum_{l\neq 0} b(l,N)^2  \;.
\end{split}
\end{equation}
Moreover, 
\begin{equation*}
b(l,N)^2 =\sum \sum_{n\cdot \Delta(x) = l = n'\cdot \Delta(x')} 
\widehat f(\frac nN)\widehat f(\frac {n'}N) \;.
\end{equation*}
Now summing over all $l\neq 0$ we get 
\begin{equation}\label{sum of b(l,N)}
\sum_{l\neq 0} b(l,N)^2 = \sum \sum_{n\cdot \Delta(x) = n'\cdot \Delta(x')} 
\widehat f(\frac nN)\widehat f(\frac {n'}N)
\end{equation}

Fix $\epsilon>0$, and choose $\delta = \epsilon/2k$ and
$R$ sufficiently large in terms of $k$ and $\delta$, say
$R>2k+(4k+100)/\delta$.   
Also set $M=N^{1+\delta}$. 
We have $\widehat f(x)\ll |x|^{-R}$ for large $x$.
In \eqref{sum of b(l,N)} we break up the sum  over $n$ into ranges
$0<|n| \leq M$ and $|n|>M$, 
and likewise for the sum over $n'$.  
In the range $0<|n|<M$ we use the bound $|\widehat f(\frac nN)|\ll 1$,
and in the range $|n|>M$ we use $\widehat f(x)\ll |x|^{-R}$. This
gives 
\begin{equation}\label{big sum}
\begin{split}
\sum_{l\neq 0} b(l,N)^2 &\ll_f \sumstar_{x_i\leq N}\sumstar_{x'_i\leq N} 
\#\{0<|n|,|n'|\leq M, n\cdot \Delta(x)=n'\cdot \Delta(x') \}  \\
&+\sumstar_{x_i\leq N}\sumstar_{x'_i\leq N} 
 \sum_{0<|n|\leq M}\sum_{|n'|>M}|\frac {n'}N|^{-R}  \\
&+ \sumstar_{x_i\leq N}\sumstar_{x'_i\leq N} 
\sum_{|n|>M}|\frac {n}N|^{-R}\sum_{|n'|>M}|\frac {n'}N|^{-R} \;.
\end{split}
\end{equation}

The third term in \eqref{big sum} is bounded by square of the number of
$x_i\leq N$ times 
the square of the sum $\sum_{|n|>M}|\frac {n}N|^{-R}$, giving a total
of at most 
$$N^{2k}N^{2R}M^{-2(R-k)} \ll N^{-100}\;.
$$

The second term in \eqref{big sum} is bounded by 
\begin{equation*}
\begin{split}
N^{2k}\#\{|n|<M\} \sum_{|n'|>M}|\frac {n}N|^{-R} &\ll
N^{2k+R}M^{k-1-R+k}\\
&\ll N^{2k+(1+\delta)(2k-1)-R\delta}\ll N^{-100} \;.
\end{split}
\end{equation*}

The first term of \eqref{big sum} is bounded by the number of solutions
of the equation  
$n\cdot \Delta(x) = n'\cdot \Delta(x')$ in variables $0<|n|,|n'|\leq
M$, $x_i\leq M$ distinct, $x'_j\leq M$ distinct. By Lemma~\ref{lem z4}, this
number is at most $M^{2k-1}\log^{2k-1} M\ll N^{2k-1 +\epsilon}$. 

Thus we find that 
$$\sum_{l\neq 0} b(l,N)^2 \ll N^{2k-1+\epsilon}$$
and inserting into \eqref{var formula} we get
$$
\var(R_k(f,N))\ll N^{-1+\epsilon}  \;.
$$
\end{proof}

\section{Small fractional parts}\label{sec:small}

Our next goal will be almost-everywhere convergence. Preliminary to
that, we have to investigate the frequency of occurrence of fractional
parts of $\alpha a(x)$ in short (of size $1/N$) intervals. 
We denote by $||x||$ the distance to the nearest integer. Our
principal result in this section is: 
\begin{prop} \label{small fractional parts} 
Let $a(x)$ be lacunary and let $c > 1$ be such that 
$$
a(x+1) > c a(x)
$$ 
for all $x$. Then for almost all $\alpha$ the following holds true:
For any $\epsilon > 0$ there exists a constant 
$C$  depending only on $c$, $\alpha$ and $\epsilon$ such that for
any positive integer $N$ and any real number $\beta$ one has:
$$
\# \{x<N : ||\alpha a(x) -\beta || < 1/N\} < C N^{\epsilon}\;.
$$
\end{prop}

We first prove the following :
\begin{lem}\label{lem 3.2}
Let $N>1$ and $a_1,\cdots, a_k$ positive integers such that $a_{j+1}
\ge Na_j$ for $1\le j\le k-1 $. Then the set 
$$
\Lambda (\vec a,N)=\{\alpha \in [0,1];||\alpha a_j||\le \frac 1N,1\le j\le k\}
$$ 
has Lebesgue measure $\le \frac {4^k}{N^k}$.
\end{lem}
\begin{proof} 
Let  $\alpha \in \Lambda (\vec a,N)$. For $1\le j\le k$ we write 
$\alpha $ in the form 
$$
\alpha =\frac {b_j}{a_j}+\beta _j
$$
with $b_j=b_j(\alpha )\in \{0,1,\cdots, a_j\}$ and $\beta _j\le \frac {1}{2a_j}$.
 From $\alpha a_j=b_j+a_j\beta _j$, with $b_j\in {\bf Z}$ 
and $a_jb_j\in [-\frac 12,
\frac 12]$ it follows that $||\alpha a_j||=|a_j\beta _j|$ and since 
$||\alpha a_j||\le \frac 1N$ we get 
$|\beta _j|\le \frac {1}{Na_j}$ for $1\le j\le k.$
For any $j\in \{1,\cdots, k\}$ let 
$$
B_j=\{0\le b\le a_j: \, \,  \text {there is} \, \, \alpha \in 
\Lambda (\vec a,N) \, \, \text {with} \, \, b_j(\alpha )=b\}.
$$
Then for any $j$ 
$$
 \Lambda (\vec a,N) \subseteq \bigcup_{b\in B_j}\Big [\frac {b}{a_j}-
\frac {1}{a_jN},\frac {b}{a_j}+\frac {1}{a_jN}\Big ]=A_j, \, \, \text { say}.
$$
In particular one has:
$$
\meas (\Lambda(\vec a,N))\le\meas (A_k)=\frac {2}{a_kN}\# B_k.
$$
It remains to bound $\# B_k$. In order to do this we produce for any
$j$ an upper bound for $\# B_j$ in terms of $\# B_{j-1}$. Let $b\in B_j$.
There is $\alpha $ such that $b_j(\alpha )=b$.
Write:
$$
\alpha =\frac {b}{a_j}+\beta _j=\frac {b_{j-1}}{a_{j-1}}+\beta _{j-1}.
$$

Then one has :
$$
|b-\frac {a_jb_{j-1}}{a_{j-1}}|=a_j|\beta _{j-1}-\beta _j|\le a_j
(\frac {1}{Na_{j-1}}+\frac {1}{Na_j})=\frac {a_j}{Na_{j-1}}+\frac 1N.
$$
For a fixed value of $b_{j-1}$ the integer $b$ may vary in the above
interval of length $\frac {2a_j}{Na_{j-1}}+\frac 2N$, so it takes at most 
$1+[\frac {2a_j}{Na_{j-1}}+\frac 2N]\le 2+\frac {2a_j}{Na_{j-1}}$
values. Hence:
$$
\# B_j\le 2\Big (1+\frac {a_j}{Na_{j-1}}\Big )\# B_{j-1},2\le j\le k.
$$
Clearly $\# B_1\le (1+a_1)$. By multiplying these inequalities we obtain:
$$
\# B_k\le 2^{k-1}(1+a_1)(1+\frac {a_2}{Na_1})\cdots (1+\frac {a_k}{Na_{k-1}})
$$
and therefore 
$$
\meas (\Lambda(\vec a,N))\le 
\frac {2^k}{a_kN}(1+a_1)(1+\frac {a_2}{Na_1})\cdots 
(1+\frac {a_k}{Na_{k-1}})
$$
$$
=\frac {2^k}{N^k}\frac {1+a_1}{a_1}
\frac {Na_1+a_2}{a_2}\cdots \frac {Na_{k-1}+a_k}{a_k}.
$$
Here we use the  assumption that $Na_j\le a_{j+1}$ to conclude that 
$$
\meas \Lambda (\vec a,N)\le \frac {4^k}{N^k}
$$
which completes the proof of the lemma.
\end{proof}

We now introduce some notation.
Given $N\ge 1$ and $\alpha ,\beta \in [0,1]$ denote 
$$
G(N,\alpha ,\beta )=\#\{x\le N;||\alpha a(x)-\beta ||<\frac 1N\}.
$$
Then set:
$$
G(N,\alpha )=\max_{\beta \in [0,1]}G(N,\alpha ,\beta ).
$$
Given $\delta >0$ and $N\ge 1$ define the set :
$$
A(\delta ,N)=\{\alpha \in [0,1]:G(N,\alpha )>N^{\delta }\}.
$$
Note that by the above definitions, if $\alpha $ is not in the exceptional set 
$A(\delta ,N)$ then $G(N,\alpha )\le N^{\delta }$ so uniformly for all $\beta $
one has 
$G(N,\alpha ,\beta )\le N^{\delta }$, i.e.
$$
\# \{x\le N;||\alpha a(x)-\beta ||<\frac 1N\} < N^\delta
$$
for all $\beta $. Set $\tilde A(\delta ,M)=\bigcup _{N\ge M}A(\delta ,N)$ and 
$\tilde A(\delta )=\bigcap _{M\ge 1} \tilde A(\delta ,M)$.
Now let $\alpha \notin \tilde A(\delta )$.
Then there exists $M=M(\alpha ,\delta ) $ such that $\alpha $ is not in $\tilde 
A(\delta ,M)$. Thus for any $N\ge M(\alpha ,\delta )$ we have $\alpha \notin A
(\delta ,N)$ and so:
For any $N\ge M(\alpha ,\delta )$ we have uniformly for all $\beta$ :
$$
\# \{x\le N;||\alpha a(x)-\beta ||<\frac 1N\}\le N^{\delta }.
$$
In other words, if $\alpha \notin \tilde A(\delta )$ then there exists 
$C(\delta ,\alpha ,c)$ such that for all $N$ and all $\beta $ one has :
$$
\# \{x\le N;||\alpha a(x)-\beta ||<\frac 1N\}\le C(\delta ,\alpha ,c)N^{\delta }.
$$
In order to prove Proposition~\ref{small fractional parts} 
we need to show that for any $\delta >0$
the set $\tilde A(\delta )$ has measure zero. Fix $\delta >0$.
By the definition of $\tilde A(\delta )$ one has 
$\meas \tilde A(\delta )\le \meas \tilde A(\delta ,M)$ for any $M\ge 1$, 
so it is enough to show that:
\begin{equation}\label{z1}
\meas \tilde A(\delta ,M)\rightarrow 0 \, \, \text {as}\,\,
 M\rightarrow \infty .
\end{equation}
Now $\meas  \tilde A(\delta ,M)\le \sum _{N\ge M} \meas A(\delta ,N)$.
Thus in order to prove (\ref{z1}) it is enough to show that 
there exists $\epsilon _{\delta }>0$ such that for any $N\ge 1$ one has:
\begin{equation}\label{z2}
\meas A(\delta ,N)\ll_{c,\delta }\frac {1}{N^{1+\epsilon_ \delta }}.
\end{equation}
We will prove this in the next Lemma, 
which completes the proof of Proposition~\ref{small fractional parts}.  
\begin{lem}
Given $\delta >0$, for any $N\ge 1$ one has :
$$
\meas A(\delta ,N)\ll_{c,\delta }\frac {1}{N^{1999}}.
$$
\end{lem}
\begin{proof} 
Given $\delta >0$ we choose a positive integer $k$,  depending on $\delta $
only, whose precise value will be given later. Let $N\ge 1$ and $\alpha 
\in A(\delta ,N)$.
There exists $\beta \in [0,1]$ such that the set 
$$
\mathcal N=\{x\le N;||\alpha a(x)-\beta ||<\frac 1N\}
$$
has more than $[N^{\delta }]$ elements. Arrange the elements of $\mathcal N$
in increasing order: $\{1\le x_1<x_2<\cdots <x_l\}$ and pick from this set 
the first element $x_1$, then ignore the next $r=[(1+\delta )\log _cN]$
elements, pick the next one, ignore again $r$ elements, and so on.
We get a set of ``well spaced'' integers 
$\mathcal M=\{y_1=x_1<y_2=x_{r+1}<y_3<\cdots < y_s\}$
with $s\ge \frac {N^{\delta }}{1+(1+\delta )\log _cN}\, $, 
such that 
$$
||\alpha a(y_j)-\beta ||<\frac 1N, \qquad 1\le j\le s
$$
and (since $y_{j+1}-y_j\geq (1+\delta) \log_c N $):  
\begin{equation}\label{z3}
a(y_{j+1})\ge N^{1+\delta }a(y_j),\qquad 1\le j\le s-1 \;.
\end{equation}
Now look at the sequence of fractional parts $\mathcal U=
(\{\alpha a(y_j)\})_{1\le j\le s}$.
They all fall in an interval of length $\frac 2N$ centered in $\{\beta \}$.
We cut this interval in $m=[\frac {s-1}{k}]$ intervals $J_1,\cdots ,J_m$
having the same length: $\frac {2}{Nm}$. By the box principle, one of these 
intervals, $J_{i_0}$ say, will contain at least 
$\frac sm=\frac {s}{[\frac {s-1}{k}]}>\frac {s-1}{[\frac {s-1}{k}]}\ge k$
elements of $\mathcal U$, that is, $J_{i_0}$ will contain at least $k+1$ 
elements of $\mathcal U$. So let $z_0<z_1<\dots <z_k$ be 
$k+1$ elements of $\mathcal M$ for which the fractional parts 
$\{\alpha a(z_0)\},\cdots ,\{\alpha a(z_k)\}$ belong to $J_{i_0}$.
Then clearly one has:
\begin{multline}\label{z4}
||\alpha (a(z_1)-a(z_0))||,\dots ,
||\alpha (a(z_k)-a(z_0))|| \le \text {length} |J_{i_0}|\\
=\frac {2}{Nm}=\frac {2}{N[\frac {s-1}{k}]}\le \frac {4k}{Ns}
\le \frac {4k(1+(1+\delta )\log _cN)}{N^{1+\delta }}<
\frac {1}{N^{1+\frac {\delta }{2}}}
\end{multline}
for $N$ sufficiently large  in terms of $c,k$ and $\delta $.
Note also that since the $z_i$ are still well-spaced, 
by (\ref{z3}) one has:
\begin{equation}\label{z5} 
a(z_1)\ge N^{1+\delta }a(z_0),\dots ,a(z_k)\ge N^{1+\delta }
a(z_{k-1}) \;.
\end{equation}
Let $\vec a =(a_1,\dots ,a_k)$ be given by:
$$
a_1=a(z_1)-a(z_0),\dots ,a_k=a(z_k)-a(z_0) \;.
$$
By \eqref{z5} we see that for $i=1,\cdots ,k-1$ one has:
\begin{multline}\label{z6}
a_{i+1}=a(z_{i+1})-a(z_0)\ge N^{1+\delta}a(z_i)-a(z_0) \\ 
> N^{1+\delta}(a(z_i)-a(z_0))
=N^{1+\delta }a_i
\end{multline}
while (\ref{z4}) says that 
\begin{equation}\label{z7}
||\alpha a_i||<\frac {1}{N^{1+\frac {\delta }{2}}},1\le i\le k.
\end{equation}
 From (\ref{z6}) and (\ref{z7}) we see that one may apply 
Lemma~\ref{lem 3.2} to the vector  
$\vec a$, with $N$ replaced by $N^{1+\frac {\delta }{2}}$. In the terminology of 
that Lemma, $\alpha $ belongs to $\Lambda (\vec a,N^{1+\frac {\delta }{2}})$.
Since for each $\alpha \in A(\delta ,N)$ there is such a vector $\vec a $ it 
follows that 
$$
A(\delta ,N)\subseteq \bigcup _{\vec a}\Lambda (\vec a,N^{1+\frac {\delta }{2}}).
$$ 
By Lemma~\ref{lem 3.2} we derive:
$$
\meas A(\delta ,N)\le \sum _{\vec a}
\meas\Lambda (\vec a,N^{1+\frac {\delta }{2}}) 
\le \frac {4^k\#\{\vec a\}}{N^{k(1+\frac {\delta }{2})}}.
$$
Now each vector $\vec a $ as above is uniquely determined by a $(k+1)$-tuple 
$(z_0,z_1,\dots, z_k)$ of positive integers $\le N$. 
The number of such $(k+1)$-tuples is $<N^{k+1}$.
It follows that
$$
\meas A(\delta ,N)\ll_{c,k,\delta }
\frac {N^{k+1}}{N^{k(1+\frac {\delta }{2})}}=
\frac {N}{N^{\frac {k\delta }{2}}}.
$$
We now let $k=\frac {4000}{\delta }$ and the lemma is proved.
\end{proof}

\section{Almost everywhere convergence} \label{sec:ae}

We now show that there is a set of  $\alpha$ of full measure so that for 
all $k\geq 2$ and all test functions $f\in C_c^\infty(\R^{k-1})$, the
$k$-level correlation functions $R_k(f,N)(\alpha)$ converge to 
$\int f(x)dx$.  
The main ingredient here is:
\begin{prop}\label{prop X}
Fix  $f\in C_c^\infty(\R^{k-1})$. 
If  $0<\delta<1$ and $1\leq K\leq N^{1-\delta}$ then for almost every $\alpha$ 
$$
R_k(f,N+K)(\alpha)-R_k(f,N)(\alpha) \to 0 
$$
\end{prop}

\subsection{Proof of Theorem~\ref{ae for R_k}}
We first show how Proposition~\ref{prop X} implies 
Theorem~\ref{ae for R_k}: 
By Proposition \ref{prop:variance}, for fixed $f$ we have  
$$
\int_0^1 \left| R_k(f,N)(\alpha) -\EX(R_k(f,N)) \right|^2 d\alpha \ll_\epsilon
N^{-99/100 }
$$ 
and so if we take  $N_m \sim m^{101/99}$ then 
\begin{multline*}
\int_0^1 \sum_m \left| R_k(f,N_m)(\alpha) -\EX(R_k(f,N) \right|^2
d\alpha \\
= \sum_m \int_0^1 \left| R_k(f,N_m)(\alpha) -\EX(R_k(f,N) \right|^2 d\alpha \\
<\sum_m \frac 1{m^{101/100}}  <\infty
\end{multline*}
Thus the sum $\sum_m |R_k(f,N_m)(\alpha) -\EX(R_k(f,N)|^2$ is finite almost
everywhere,  and hence the individual summands converge to zero as
$m\to\infty$ for almost all $\alpha$.  

For each $N$ we can find $m$ such that $N_m\leq N <N_{m+1}$.  
Then since $R_k(f,N_m)(\alpha)- \EX(R_k(f,N)\to 0$ for almost all
$\alpha$, and by Lemma~\ref{lem:average}, $\EX(R_k(f,N) \to \widehat f(0)$, 
Proposition~\ref{prop X} will
show that $R_k(f,N)(\alpha)\to\widehat f(0) $ for a set of full
measure of $\alpha$ which depend on the test function $f$.
By a standard diagonalization argument one can pass to a subset of
full measure of $\alpha$'s which work for all $f$'s (see \cite{RS}). \qed

\subsection{An upper bound for $R_k(f,N)$} 
As a consequence of Proposition~\ref{small fractional parts} we have the 
following a-priori estimate on the correlation functions: 
\begin{lem}\label{few small points}
 For almost all $\alpha$ we have 
$$
R_k(f,N)(\alpha)\ll_{\epsilon,f} N^\epsilon
$$
\end{lem}
\begin{proof}
We use the representation of $R_k(f,N)$ as in
\eqref{definition of R_k}:  
$$
R_k(f,N)(\alpha):=\frac 1N\sumstar_{x_i\leq N} F_N(\alpha \Delta(x))
$$
where $\Delta(x) := (a(x_1)-a(x_2),\dots, a(x_{k-1})-a(x_k))$.  
Note that 
$$|R_k(f,N)|\leq R_k(|f|,N)
$$ 
so we may assume $f\geq 0$. Now fix $x_1$, and set $\beta=\alpha a(x_1)$; 
then for $\alpha \Delta(x)$ to lie in the support of $F_N$, we need 
$||\alpha a(x_2) -\beta|| \ll_f 1/N$.  
By Proposition~\ref{small fractional parts}, for almost all $\alpha$  
there are at most $O_{f,\epsilon}(N^\epsilon)$ integers $x_2\leq N$
satisfying this. Similarly, we need $||\alpha a(x_i) -\beta|| \ll_f 1/N$ 
for all $2\leq i\leq k$ which forces the number of possible
$x=(x_1,\dots x_k)$ 
contributing to the sum to be at most $O(N^\epsilon)$. 
Now summing over the $N$ possible $x_1$'s gives 
$R_k(f,N)(\alpha)\ll_{\epsilon,f} N^\epsilon$. 
\end{proof}

\subsection{Proof of Proposition~\ref{prop X}} 
Now fix $0<\delta<1$ and assume that $K\leq N^{1-\delta}$.
We will show that for almost all $\alpha$, 
$$
\left| R_k(f,N+K)(\alpha)-R_k(f,N)(\alpha) \right|  \ll KN^{-1+\epsilon}
$$

\noindent{\bf Step 1:} {\em In the expression }
$$
R_k(f,N+K)=\frac 1{N+K}\sumstar_{x_i\leq N+K} F_{N+K}(\alpha \Delta(x))
$$
{\em we can replace $1/(N+K)$ by $1/N$ with error
$O_{\epsilon,f}(KN^{-1+\epsilon})$. }
 
Indeed,by Lemma~\ref{few small points}, 
$R_k(f,N+K)\ll N^\epsilon$ and so 
\begin{equation*}
\begin{split}
\frac 1N\sumstar_{x_i\leq N+K} F_{N+K}(\alpha \Delta(x)) &=
(1+\frac KN) R_k(f,N+K)\\  
&= R_k(f,N+K)+O( \frac KN N^\epsilon)
\end{split}
\end{equation*}
as claimed.  

\noindent{\bf Step 2:} {\em We may replace the sum over (distinct)
$x_i\leq N+K$ by the sum over (distinct) $x_i\leq N$: }
$$
\sumstar_{x_i\leq N+K} F_{N+K}(\alpha \Delta(x)) = 
\sumstar_{x_i\leq N} F_{N+K}(\alpha \Delta(x)) + O(KN^\epsilon) \;.
$$
Indeed, the difference between the two sums is a sum over a union of
subsets 
$$
S(I) = \{(x_1,\dots,x_k)\text{ distinct}: N<x_i\leq N+K ,  i\in
I, x_j\leq N,  j\notin I \}
$$
where the index set $I$ runs over all the $2^{k}-1$ nonempty subsets
of $\{1,2,\dots,k\}$. 

To estimate the contribution of 
$\Sigma(I):=\sum_{x\in S(I)}F_{N+K}(\alpha\Delta(x)) $, we use the
consequence of Proposition~\ref{small fractional parts}, which says that if
we fix one of the coordinate axes $i_0$, then the number of vectors
$x$ with $x_{i_0}=y$  
fixed which contribute to the sum is $O(N^\epsilon)$, uniformly in
$y$. Thus the number of vectors in $S(I)$ which contribute to the sum
$\Sigma(I)$ is at most $O(KN^\epsilon)$, because if we look at $i_0\in
I$ we have $N<x_{i_0}\leq N+K$ for $x\in S(I)$, and so 
$$
\Sigma(I) \ll_{f,\epsilon} KN^\epsilon \max|f| \ll KN^\epsilon \;.
$$
Thus we find 
\begin{equation*}\begin{split}
R_k(f,N+K)-R_k(f,N)  &= \frac 1N\sumstar_{x_i\leq N} 
F_{N+K}(\alpha \Delta(x))  - F_N(\alpha \Delta(x))\\ 
& + O_{f,\epsilon}(KN^{-1+\epsilon} ) \;.
\end{split}
\end{equation*}

\noindent{\bf Step 3:} 
{\em  We show that for almost all $\alpha$,}
$$
\frac 1N\sumstar_{x_i\leq N} 
F_{N+K}(\alpha \Delta(x))  - F_N(\alpha \Delta(x)) \ll
KN^{-1+\epsilon}\;.
$$

{\noindent}{\em Remark:} This is the statement that the correlation
functions are independent of the exact unfolding procedure! 

First, a digression: Given a vector $y\in \R^{k-1}$, there is a unique
integer vector $m_y\in \Z^{k-1}$ so that $y+m_y$ lies in the cube
$(-1/2,1/2]^{k-1}$. Moreover, for any other integer vector  $m\neq
m_y$, $||m+y||>1/2$ and so $||N(m+y)||>N/2$.  
Thus if $N$ is sufficiently large so that $\supp(f)$ lies in a ball of
radius $\rho(f)<N/2$ around the origin, then 
$$
F_N(y) = f(N(m_y+y))
$$
and 
\begin{equation}\label{ball}
||N(m_y+y)||<\rho(f) \;.
\end{equation}
Furthermore, if $m\neq m_y$ then $||(N+K)(m_y+y)||>||N(m+y)||>N/2$ and
therefore
$$
F_{N+K}(y) = f((N+K)(m_y+y)) \;.
$$

Apply these considerations to $y=\alpha \Delta(x)$ and abbreviate
$$
v_x:=m_{\alpha \Delta(x)}+\alpha\Delta(x)
$$ 
to get that if $N>N_0(f)$ then 
\begin{multline*}
\frac 1N\sumstar_{x_i\leq N} 
F_{N+K}(\alpha \Delta(x))  - F_N(\alpha \Delta(x)) \\= 
\frac 1N\sumstar_{x_i\leq N} 
f((N+K)v_x) - f(N v_x) \;.
\end{multline*}
By the mean value theorem, 
\begin{equation}\label{mean value} 
\begin{split}
f((N+K)v_x) - f(N v_x)
&= Kv_x \cdot \nabla f(N v_x +\theta K v_x )\\
&= \frac KN N v_x\cdot  
\nabla f(N v_x (1+\theta \frac KN))
\end{split}
\end{equation}
for some $0<\theta=\theta_x <1$ depending on $x$. 
If this is nonzero, then certainly $Nv_x$ is contained in a ball of
radius $2\rho(f)$ around the origin. 
Now $||Nv_x||<\rho(f)$ by \eqref{ball}, so the sum of the terms
\eqref{mean value}  is bounded by $\rho(f)\max ||\nabla f||$ times the
number of $x$ for which $Nv_x$ lies in a ball of radius $2\rho(f)$
around the origin. 

We can now bound the sum of \eqref{mean value} 
by relating it to a
smoothed $k$-level correlation function as follows: 
Choose a positive, smooth function $g\in C_c^\infty(\R^{k-1})$ which
is constant on the ball of radius $2\rho(f)$ around the origin, and
satisfies $g\geq \max ||\nabla f||$.  Write $G_N(y):=\sum_m g(N(m+y))$. 
Then 
$$
\nabla f(N v_x (1+\theta \frac KN))\leq 
g(N v_x) = G_N(\alpha \Delta(x)) \;.
$$
Thus we find that 
\begin{equation*}\begin{split}
\frac 1N\sumstar_{x_i\leq N} 
f((N+K) v_x) - f(N v_x)
&\ll \frac KN \rho(f) 
\frac 1N\sumstar_{x_i\leq N}G_N(\alpha \Delta(x)) \\
&=\frac KN \rho(f) R_k(g,N) \;.
\end{split}\end{equation*}
By Lemma~\ref{few small points}, $R_k(g,N)\ll_{g,\epsilon} N^\epsilon$
for a.e.  $\alpha$, which gives the result of step $3$. 
This concludes the Proof of Proposition~\ref{prop X}. \qed

\end{document}